\documentclass{article}
\usepackage{arxiv}
\usepackage{graphicx} 
\usepackage{amssymb}
\usepackage{color}
\usepackage[T1]{fontenc}
\usepackage{amsthm}
\usepackage{amsmath,amsfonts}
\usepackage{enumitem}
\usepackage{cite}
\usepackage{hyperref}

\newtheorem{definition}{Definition}
\newtheorem{proposition}{Proposition}
\newtheorem{theorem}{Theorem}
\newtheorem{remark}{Remark}

\title{Pattern formation: reactivity is not necessary for chemotaxis--driven instabilities}

\author{
 Angela Monti \\
  Istituto per le Applicazioni del Calcolo \lq \lq M. Picone\rq \rq \\
  National Research Council (CNR)\\
  via G. Amendola 122/D, Bari, Italy\\
  \texttt{angela.monti@cnr.it} }

\begin{document}

\maketitle
\begin{abstract}
A classical result by Neubert, Caswell and Murray states that reactivity of a spatially homogeneous equilibrium is a necessary condition for diffusion-driven (Turing) instability. In this work, we investigate whether the same conclusion remains valid in the presence of chemotaxis. We consider a general reaction--diffusion system coupled with a chemotactic flux and establish necessary conditions for asymptotic instability. We show that the classical requirement of reactivity can be relaxed when the chemotactic contribution is sufficiently strong. In particular, while reactivity remains necessary for Turing instability, it is not a necessary condition for chemotaxis-driven instability.
From a computational viewpoint, we extend the matrix-oriented formulation developed for reaction--diffusion systems to the more general class of reaction--diffusion--chemotaxis models. The chemotactic transport term is discretized in a form compatible with the matrix-oriented approximation of the diffusion operator, yielding an efficient numerical framework for the simulation of chemotaxis-driven pattern formation.
A geometric interpretation of the instability region is presented, highlighting the distinct roles played by diffusion and chemotaxis. The theoretical and numerical developments are illustrated through two representative examples: a chemotaxis-extended Schnakenberg model, showing how chemotaxis modifies classical Turing patterns, and a predator--prey model, demonstrating that chemotaxis alone can induce pattern formation in the absence of both reactivity and diffusion-driven instability.
These results reveal a fundamental difference between diffusion-driven and chemotaxis-driven mechanisms of spatial self-organization and provide new theoretical and computational insights into the role of non-symmetric transport processes in biological pattern formation.
\end{abstract}

\section{Introduction} 
The emergence of spatially heterogeneous structures from initially homogeneous states is one of the central topics in mathematical biology, ecology, chemistry, and the study of complex systems. Since the seminal work of Turing~\cite{Turing1952}, reaction--diffusion mechanisms have provided one of the most influential explanations for spontaneous pattern formation through the loss of stability of a homogeneous equilibrium. In this framework, a spatially homogeneous steady state which is stable for the local reaction kinetics may become unstable in the presence of diffusion, leading to the formation of stationary spatial patterns. This mechanism has been extensively studied in chemical and biological systems~\cite{Murray_book2,Maini1997, CrossHohenberg1993}, and it underlies a wide range of models for morphogenesis, ecological segregation, and self-organization phenomena. Although the classical Turing mechanism remains a central paradigm, recent developments have emphasized that pattern formation cannot always be fully understood within the standard two-species reaction--diffusion setting. Extensions of Turing's theory include growing and evolving domains, heterogeneous media, localized structures, cross-diffusion, nonlocal interactions, and more general reaction--transport systems ~\cite{madzvamuse2010stability,AlSaadi2021,Vanag2007,Brena2014,Gambino2012Turing,Gambino2013CrossDiffusion, VillarSepulveda2025ReactionCrossDiffusion,Ninomiya2017Nonlocal}. Moreover, it has recently been shown that Turing-type linear instability conditions alone do not necessarily guarantee the emergence of persistent nonlinear patterns~\cite{Krause2023}. This motivates the combined use of linear stability analysis, transient dynamics, and numerical simulations in the study of spatial self-organization. A related but distinct mechanism arises in systems where movement is not purely diffusive. In many biological and ecological contexts, individuals, organisms, or cells move directionally in response to environmental cues. Chemotaxis, namely the biased movement induced by chemical or biological gradients, is a classical example of such directed transport, with roots in the works of Patlak and Keller--Segel~\cite{Patlak1953, keller1970initiation,KellerSegel1971}. Chemotaxis and related taxis mechanisms have been widely used to model bacterial aggregation, cellular migration, predator--prey interactions, tumor invasion, and ecological pattern formation ~\cite{Hillen2008,Bellomo_2022}. In contrast with ordinary diffusion, chemotaxis introduces a non-symmetric transport contribution, often related to cross-diffusion effects, whose destabilizing or stabilizing role depends on the coupling between taxis sensitivity and the local reaction kinetics ~\cite{CaoWu2021}. 

The starting point of the present work is the relationship between pattern formation and reactivity. For an asymptotically stable equilibrium, reactivity measures the possibility of short-time amplification of perturbations before their eventual asymptotic decay. In finite-dimensional systems, this quantity is given by the numerical abscissa, i.e., the largest eigenvalue of the Hermitian part of the linearized operator ~\cite{NeubertCaswell1997,trefethen2005spectra}. Although Turing instability is an asymptotic mechanism and reactivity is a transient property, Neubert, Caswell and Murray~\cite{Neubert2002} established a remarkable connection between the two concepts: for classical reaction--diffusion systems, reactivity of the spatially homogeneous equilibrium is a necessary condition for diffusion-driven, or Turing, instability. In other words, diffusion alone cannot destabilize a stable non-reactive equilibrium. The presence of chemotaxis substantially modifies this scenario. Unlike diffusion, the chemotactic contribution does not simply add a symmetric dissipative operator to the linearized dynamics. Instead, it introduces a non-symmetric transport term whose effect depends on the interaction between the chemotactic sensitivity, the equilibrium density, and the off-diagonal structure of the kinetic Jacobian. This naturally raises the following question: 
\begin{center} 
\emph{Is reactivity still a necessary condition for pattern formation when chemotaxis is present?} 
\end{center} 

The main purpose of this work is to answer this question for a general two-species reaction--diffusion--chemotaxis model. Building on the well known linear instability conditions for reaction--diffusion--chemotaxis systems~\cite{CaoWu2021}, we analyze whether the reactivity requirement established by Neubert et al.~\cite{Neubert2002} for classical Turing instability remains valid when chemotactic transport is present. This comparison reveals that the Neubert--Caswell--Murray result does not extend in a direct way to chemotactic systems. More precisely, while reactivity of the homogeneous equilibrium remains necessary for classical diffusion-driven instability, it is not necessary for chemotaxis-driven instability. Sufficiently strong chemotactic transport can destabilize a homogeneous equilibrium even when the local kinetic Jacobian is stable and non-reactive. We emphasize that the present question concerns the reactivity of the local kinetic Jacobian, corresponding to the homogeneous mode. This should be distinguished from the reactivity of the spatially dependent operators, or of the fully discretized spatial system, which is the object of several studies on transient amplification and non-normality in spatially extended models~\cite{Klika2017, muolo2019patterns,DieleBMAB2025}. 

A second contribution of the paper is computational. Numerical simulation of pattern-forming systems often requires fine spatial grids and long-time integration, especially when stationary patterns are selected after long transients. For classical reaction--diffusion systems, matrix-oriented methods exploit the tensor-product structure of the discrete diffusion operator on rectangular domains, avoiding the explicit construction of the large vectorized spatial operator. In particular, the approach introduced in ~\cite{DautiliaSguraSimoncini2020} reformulates implicit diffusion steps in terms of Sylvester matrix equations, leading to efficient matrix-based time integrators for reaction--diffusion PDEs. Other tensor-oriented strategies have also been developed for stiff multidimensional advection--diffusion--reaction systems. For instance, in ~\cite{CALIARI2024112640}, the authors proposed a second-order directional split exponential integrator based on the directional splitting of matrix functions and tensor-matrix products, with applications to pattern-forming systems. In the present work, we extend the matrix-oriented formulation, proposed in \cite{DautiliaSguraSimoncini2020}, to reaction--diffusion--chemotaxis systems. The main difficulty is the discretization of the nonlinear chemotactic flux, which must be incorporated without destroying the matrix-oriented structure of the implicit diffusion step. We introduce a compatible discretization of the chemotactic operator which reduces to the same discrete Laplacian used for diffusion in the constant-mobility limit, in the spirit of the Extended Central Difference Formulas (ECDF, \cite{amodio2005}) discretizations used for Turing pattern approximation ~\cite{SGURA20124132}. 

The theoretical and numerical developments are illustrated through two classes of examples. First, we consider a chemotaxis-extended Schnakenberg model~\cite{Schnakenberg1979}. In this case, the underlying reaction--diffusion system already exhibits classical Turing instability. Chemotaxis therefore acts by modifying, enhancing, weakening, or reshaping pre-existing Turing patterns. The parameter regimes used in the numerical tests are inspired by the interactive simulations provided in VisualPDE~\cite{walker2023visualpde}. The experiments show how the sign and magnitude of the chemotactic sensitivity affect both spot-like and stripe-like patterns, producing transitions between different spatial morphologies. Second, we consider a predator--prey model with prey-taxis, inspired by recent studies of taxis-driven spatial patterns in predator--prey systems~\cite{Wang2021,Wang2026}. This example addresses the main theoretical message of the paper. We select parameter regimes in which the homogeneous equilibrium is stable, non-reactive, and not destabilized by diffusion alone. Consequently, any spatial structure observed in the numerical simulations must be attributed to the chemotactic mechanism. This example demonstrates that chemotaxis may generate patterns in the absence of both reactivity and classical Turing instability. 

The paper is organized as follows. Section~\ref{sec:chem_insta} derives the linear instability conditions for diffusion--chemotaxis systems and establishes the corresponding necessary conditions based on reactivity arguments. Section~\ref{sec:notnecessary} explains why reactivity of the homogeneous equilibrium is not necessary for chemotaxis-driven instability. Section~\ref{sec:matrix_chemotaxis} presents the matrix-oriented discretization of the chemotactic term and its compatibility with the ECDF Laplacian used in the implicit diffusion step. Section~\ref{sec:numerical_examples} illustrates the theoretical results and the numerical method through the Schnakenberg and predator--prey examples. Finally, Section~\ref{sec:conclusions} summarizes the main conclusions and outlines possible future developments.

\section{Diffusion-chemotaxis driven instability}
\label{sec:chem_insta}
\noindent
We consider the general reaction-diffusion model with chemotaxis term 
\begin{equation}
\label{eq:general_hu}     
\begin{aligned}         
\partial_t u &=  D_u \Delta u  -\beta \, \nabla \cdot ( u \nabla v) + f(u,v) \\ \partial_t v &= D_v \Delta v + g(u,v) 
    \end{aligned}
\end{equation}
equipped with suitable initial and boundary conditions. Here, $D_u,D_v>0$ denote the diffusion coefficients of the two species, while $\beta \in \mathbb{R}$ measures the strength and direction of the chemotactic response.
The term $-\beta \nabla \cdot(u \nabla v)$ describes a directed movement of the species $u$ induced by spatial gradients of $v$. In conservative form, the corresponding chemotactic flux is proportional to $\beta u\nabla v$. Thus, when $\beta>0$, the population $u$ moves preferentially towards regions of high $v$-density (concentration). Conversely, when $\beta<0$, the movement is directed away from such regions.

In this section, we briefly recall the well-known linear stability analysis for diffusion-chemotaxis driven instabilities, both for completeness and to fix the notation. Such analyses appear in the literature \cite{CaoWu2021,DieleBMAB2025}, where the interplay between diffusion and chemotactic flux yields pattern formation even beyond classical Turing regimes. 

\subsection{Linear stability analysis}
We recall only the ingredients of the linear stability analysis that will be needed in the sequel. A detailed derivation of the diffusion--chemotaxis instability conditions can be found in \cite{CaoWu2021,DieleBMAB2025}. Let $(u^*,v^*)$ be a spatially homogeneous equilibrium of the reaction kinetics, namely \[ f(u^*,v^*)=g(u^*,v^*)=0, \] and denote by \begin{equation} \label{eq:jacobian} J^* = \begin{pmatrix} f_u & f_v \\ g_u & g_v \end{pmatrix}_{(u^*,v^*)} \end{equation} the kinetic Jacobian evaluated at the equilibrium. We assume that the homogeneous equilibrium is asymptotically stable in the absence of spatial transport, that is \begin{equation} \label{eq:conditions_lin} f_u+g_v<0, \qquad f_u g_v-f_v g_u>0 . \end{equation} 
Linearizing \eqref{eq:general_hu} about $(u^*,v^*)$ gives \begin{equation} \label{eq:general_full_lin} 
\mathbf{w}_t = L\Delta \mathbf{w} + J^*\mathbf{w}, \end{equation} 
where 
$$\mathbf{w}= \begin{pmatrix} 
u-u^*
\\ v-v^* 
\end{pmatrix}, \qquad 
L= \begin{pmatrix} 
D_u & -\beta u^*\\ 
0 & D_v 
\end{pmatrix}.$$
Considering Fourier modes with wavenumber $k=\|\boldsymbol\omega\|$ leads to 
\begin{equation} 
\label{eq:general_linear} 
\frac{d \widetilde{\mathbf{w}}}{dt} = J_k \widetilde{\mathbf{w}}, \end{equation} 
where 
\begin{equation} 
\label{eq:j_matrix} 
J_k = J^* - k^2L = J^* - k^2D - k^2C, 
\end{equation} 
with 
$$D = 
\begin{pmatrix} 
D_u & 0 \\ 
0 & D_v 
\end{pmatrix}, \qquad C = \begin{pmatrix} 
0 & -\beta u^* 
\\ 0 & 0 
\end{pmatrix}, \qquad L=D+C .$$ 
Thus, diffusion enters through the symmetric positive diagonal matrix $D$, whereas chemotaxis contributes through the non-symmetric matrix $C$. The eigenvalues of $J_k$ are the roots of 
\begin{equation} 
\label{eq:charpol} \lambda^2 + \lambda\left((D_u+D_v)k^2-(f_u+g_v)\right) + h(k^2) = 0, 
\end{equation} 
where 
\begin{equation} 
\label{eq:h_poly} 
h(k^2) = D_uD_vk^4 - \left(D_u g_v+D_v f_u+\beta g_u u^*\right)k^2 + \left(f_u g_v-f_v g_u\right). 
\end{equation} 
Since \eqref{eq:conditions_lin} implies that the coefficient of $\lambda$ in \eqref{eq:charpol} is positive for every $k\neq0$, instability can occur only if $h(k^2)<0$ for some nonzero wavenumber. Equivalently, the diffusion--chemotaxis instability condition can be written as 
\begin{equation} 
\label{eq:suff_conditions_cd} 
D_u g_v + D_v f_u + \beta g_u u^* > 2\sqrt{D_uD_v\left(f_u g_v-f_v g_u\right)} . \end{equation} 
Together with \eqref{eq:conditions_lin}, this condition guarantees that the homogeneous equilibrium is stable with respect to the reaction kinetics but unstable with respect to spatial perturbations. For $\beta=0$, \eqref{eq:suff_conditions_cd} reduces to the classical two-species Turing condition. For $\beta\neq0$, the chemotactic contribution shifts the instability threshold through the term $\beta g_u u^*$.

The form \eqref{eq:j_matrix} makes clear that diffusion and chemotaxis affect reactivity in fundamentally different ways. Diffusion subtracts a symmetric positive semidefinite contribution from the Hermitian part of the linearized operator, whereas chemotaxis contributes through the Hermitian part of a non-symmetric transport matrix. This distinction is essential for the necessary conditions discussed next.

\subsection{Necessary conditions for asymptotic instability}
\label{sec:NCinstability}
A classical result by Neubert et al. \cite{Neubert2002} states that reactivity of the spatial homogeneous equilibrium is a necessary condition for Turing instability in reaction-diffusion systems. This follows from the definition of Turing instability: the equilibrium must be stable in the absence of diffusion but become unstable when diffusion is introduced. Before stating the formal results, we introduce the concepts and notation required to describe the relationship between asymptotic stability and reactivity. These will be essential for proving the necessary conditions for instability.

We consider the linearized system \eqref{eq:general_full_lin} and we use the Fourier transform to analyze the transformed system \eqref{eq:general_linear}. We denote by $r(J_k)$ the numerical abscissa (reactivity) \cite{trefethen2005spectra}, defined as
\begin{equation}
\label{eq:reactivity}
r(J_k):=\max_{\|\widetilde{\mathbf{w}}_0\|\neq 0} \left[\left.\left(\dfrac{1}{\| \widetilde{\mathbf{w}}\|}\dfrac{d\\\|\widetilde{\mathbf{w}}\|}{dt}\right)\right\vert_{t=0}\right]
\end{equation} 
where $\|\cdot\|$ indicates the Euclidean norm of vectors. Therefore, the numerical abscissa can also be evaluated as the spectral abscissa of the Hermitian part of matrix $J_k$, denoted by $\alpha(H(J_k))$.

\begin{definition}
We say that $J_k$ is reactive if there exist some values of $k$, such that $r(J_k)>0$\footnote{
Compared to the definitions commonly found in the literature, here, for the sake of simplicity, we transfer the adjective 'reactive' from the homogeneous equilibrium to the linear operator evaluated at equilibrium.}. 
\end{definition}
\noindent Since the numerical abscissa is bounded from below by the spectral abscissa \begin{equation}
\label{eq:spectralabscissa}
\alpha(J_k):=\max_{i=1,2} \Re (\lambda(k)_i),
\end{equation} 
we have that 
\begin{equation}
\label{eq:spectral_numerical_relation}
 \alpha(J_k)\leq r(J_k), \qquad \forall\, k.
  \end{equation}

\noindent Note that the numerical abscissa corresponds to the spectral abscissa when
$J_k$ is normal; in this case, Equation (\ref{eq:spectral_numerical_relation}) holds with equality, and the reactivity analysis provides the same information as the linear stability analysis.
 
Since instability in the presence of spatial variations requires that, for some $k\neq0$, $\alpha(J_k) >0$, the following holds: 
\begin{itemize}
    \item if $(u^*, v^*)$ is  unstable, then $J_k$ is reactive for some values of $k\neq 0$, because $\alpha(J_k)>0\Longrightarrow r(J_k)>0$;
    \item if $(u^*, v^*)$ is stable for all $k$, then $J_k$ may still be reactive for some $k$, or not reactive at all ($r(J_k)<0$ for all $k$). 
    \end{itemize}

    \begin{proposition}\label{prop:general}
        A necessary condition for diffusion-chemotaxis driven instability is that $J_k$ is reactive for some $k\neq0$.
    \end{proposition}

\noindent
A classical result \cite{Neubert2002} stems from the above proposition.
\begin{proposition}
\label{prop:turing}
        A necessary condition for Turing (diffusion-driven) instability is that  $J^*$ is reactive.
    \end{proposition}
\begin{proof}
Suppose that the homogeneous equilibrium $(u^*, v^*)$, which is stable in the absence of diffusion ($D_u = D_v = 0$), becomes destabilized by the introduction of diffusive spatial movement in the absence of the chemotaxis term ($\beta = 0$). Then, from Proposition \ref{prop:general}, the matrix $J_k = J^* - k^2 D$ is reactive, meaning that there exists at least one $k\neq 0$ such that $r(J_k) > 0$. For this $k$, let us express the Hermitian part $H(J_k)$ of the matrix $J_k$ in terms of $H(J^*)$:
\begin{equation}
\label{sym_part_J}
    H(J_k) = H(J^*) - k^2 D
\end{equation}
By Weyl’s theorem \cite{Horn_Johnson_1985}, the largest eigenvalue of the sum is less than or equal to the sum of the largest eigenvalues of each matrix, that is:
\begin{equation}
\label{eq:jstarreactive}
r(J_k) \leq r(J^*) - k^2 d_{min}
\end{equation}
where $d_{min} = \min \{D_u, D_v\} \geq 0$. Since $r(J_k) > 0$, it follows that $r(J^*) \geq k^2\,d_{min} > 0$, which proves the result.
\end{proof}

In the more general case of $\beta \neq 0$, the chemotactic contribution modifies the mode-dependent matrix $J_k$ defined in \eqref{eq:j_matrix}. Since 
\begin{equation}
\label{eq:Hjk_CD}
\begin{array}{rcl}
  H(J_k)   &= &H(J^*- k^2 C) - k^2 D \\
     & =&  H(J^*)- k^2 H(C) - k^2 D,
\end{array}
\end{equation}
the chemotactic term may either enhance or suppress reactivity, depending on the sign and size of the symmetric part $H(C)$. This yields the following necessary condition.

\begin{proposition}
\label{prop:Jreactive_gen} 
A necessary condition for diffusion-chemotaxis-driven instability is that $ J^* - k^2 C $ is reactive for some values of $k\neq0$. 
Moreover, if 
$$|\beta|u^*\le 2\sqrt{D_uD_v},$$
then reactivity of $J^*$ is also necessary for diffusion--chemotaxis-driven instability.
\end{proposition}

\begin{proof} 
Assume that the homogeneous equilibrium, stable in the absence of spatial effects, becomes unstable because of the combined action of diffusion and chemotaxis. By Proposition \ref{prop:general}, there exists a mode $k\neq0$ such that $r(J_k)>0.$ Since 
$J_k=J^*-k^2(D+C),$ its Hermitian part satisfies $H(J_k)=H(J^*-k^2C)-k^2 D$. Applying Weyl's inequality yields 
$$r(J_k) = \lambda_{\max}(H(J_k)) \le \lambda_{\max}(H(J^*-k^2C)) -k^2 d_{\min},$$ 
where $d_{\min}=\min\{D_u,D_v\}.$
Hence,
$$r(J_k) \le r(J^*-k^2C).$$ 
Since $r(J_k)>0$, it follows that $r(J^*-k^2C)>0.$ 
Therefore, $J^*-k^2C$ is reactive for some $k\neq0$. \\
To derive a condition involving $J^*$ itself, observe that $H(J_k) = H(J^*) - k^2\bigl(D+H(C)\bigr).$
Since 
$$D+H(C) = \begin{pmatrix} D_u & -\frac{\beta u^*}{2}\\ -\frac{\beta u^*}{2} & D_v \end{pmatrix},$$ 
the condition 
$$|\beta|u^* \le 2\sqrt{D_uD_v}$$ 
is equivalent to 
$$D+H(C)\succeq0,$$
that is, $D+H(C)$ is positive semidefinite. Consequently, $$H(J_k)\preceq H(J^*)$$ and, therefore, 
$$r(J_k)\le r(J^*).$$ 
Since instability requires $r(J_k)>0$ for some $k\neq0$, we obtain $r(J^*)>0.$ Hence $J^*$ is reactive. \end{proof}

\subsection{The role of chemotaxis in pattern formation}
\label{sec:chemo}
\begin{figure}
    \centering
\includegraphics[width=0.49\linewidth]{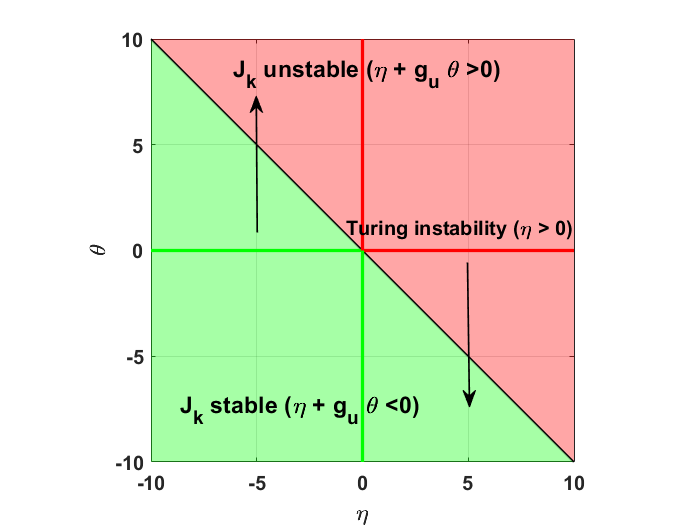}
\includegraphics[width=0.49\linewidth]{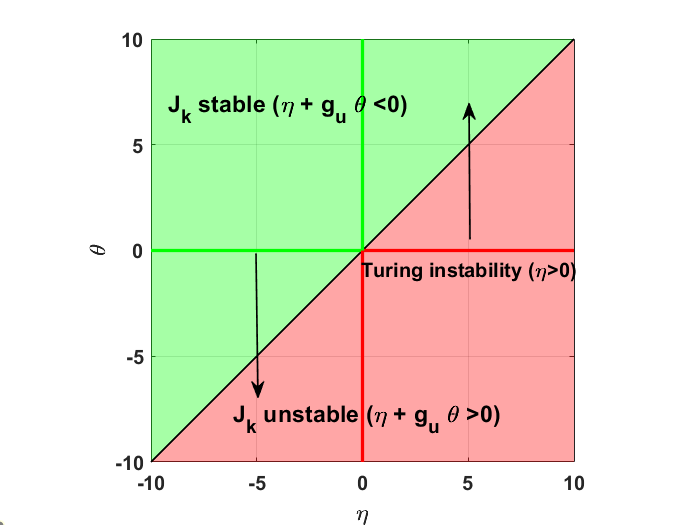}
    \caption{Stability regions for a diffusion-chemotaxis model. Left panel: case $g_u > 0$. A positive chemotaxis term $\theta = \beta u^* > 0$ can destabilize a homogeneous equilibrium that is stable in a purely diffusive regime (left arrow). Conversely, a negative chemotaxis term $\theta = \beta u^* < 0 $ can stabilize a Turing-unstable equilibrium (right arrow). Right panel: case $g_u < 0$, where the effects of positive and negative chemotaxis are reversed.
} \label{fig:primo_disegnino}
\end{figure}

To investigate the role of chemotaxis in pattern formation, we focus on the scenario associated with $J_k$-instability, i.e. \eqref{eq:suff_conditions_cd}, using the graphical representation in Figure \ref{fig:primo_disegnino}. 

Let $\eta = D_u g_v + D_v f_u - 2\sqrt{D_u D_v} \sqrt{det(J^*)}$ and $\theta = \beta \, u^*$. In the $(\eta, \theta)$ plane, the line $\eta + g_u\,\theta = 0$ separates the region where $J_k$ is unstable ($\eta + g_u\,\theta > 0$) from the region where it is stable ($\eta + g_u \theta < 0$). This representation includes, as special cases, the stability of $J^* - k^2 D $ ($\eta < 0 $) and its instability (Turing instability), represented by the half-line $\eta > 0$. Similarly, the stability of $ J^* - k^2 C $ depends on the $\theta$-axis: for $g_u > 0$, the half-line $\theta < 0$ represents stability, while $\theta > 0$ corresponds to instability; for $g_u < 0 $, these roles are reversed.

The left panel of Figure \ref{fig:primo_disegnino}
illustrates the case $g_u = 1$, highlighting how the addition of chemotaxis to a pure diffusive model can stabilize or destabilize the homogeneous equilibrium. We thus consider a point on the $\eta$-axis.

\noindent {\bf Case} $\eta>0$ (red line): when $\theta = 0$, diffusion alone destabilizes the equilibrium (Turing instability). Adding a sufficiently large negative chemotaxis sensitivity ($\theta < 0$), can restore stability (downward-pointing arrow), with the consequence of the disappearance of diffusion-driven patterns.
Conversely, adding a positive chemotaxis movement results in an enforcement of the instability as we move further away from the stable region represented by $\eta + g_u \,\theta < 0$.\\
{\bf Case} $\eta < 0$ (green line): when $\theta = 0$, diffusion is not able to destabilize the homogeneous equilibrium $(u^*, v^*)$ (for example in case when $D_u=D_v$). Adding a sufficiently large positive chemotaxis sensitivity ($\theta > 0$) may induce instability (see the upward-pointing arrow), resulting in the destabilization of the homogeneous equilibrium (chemotaxis-driven instability). In contrast, adding a negative chemotaxis movement further enhances stability, as we move further away from the unstable region represented by $\eta + g_u \,\theta > 0$. \\
Analogous considerations apply when $g_u<0$ (right panel). 

A useful remark, which will be relevant in the next sections, follows from the stability and non-reactivity of $J^*$, discussed in \cite{DieleBMAB2025} (see Proposition 1). Suppose that $J^*$ is stable and not reactive, namely
$$
f_u g_v > \frac{(f_v + g_u)^2}{4} > 0.
$$
Under this assumption, $f_u$ and $g_v$ must have the same sign. Since $J^*$ has negative trace (because of its stability), both are necessarily negative. As a consequence, the quantity $D_u g_v + D_v f_u$ is strictly negative and $\eta<0$. 

In conclusion, when the goal is to isolate the destabilizing effect produced by chemotaxis (chemotaxis‑driven instability), the analysis must be restricted to the portion of the $J_k$ instability region defined by $\eta + g_u 
\theta> 0$ that lies in the half-plane $\eta<0$. This corresponds to the regime in which the purely diffusive operator $J^*-k^2 D$ is stable, so that any loss of stability can be attributed solely to the chemotactic contribution.

\section{Reactivity is not necessary for chemotaxis-driven instability}
\label{sec:notnecessary}

A fundamental difference between diffusion-driven and
chemotaxis-driven pattern formation emerges from Proposition
\ref{prop:Jreactive_gen}.
In the purely diffusive case ($\beta=0$), Proposition
\ref{prop:turing} recovers the classical result by Neubert et
al.~\cite{Neubert2002}, stating that reactivity of the
homogeneous equilibrium is necessary for Turing instability.
Consequently, diffusion alone cannot destabilize a stable
non-reactive equilibrium.

The situation changes when chemotaxis is introduced. Indeed, the
operator governing the linearized dynamics becomes $J_k$, given by \eqref{eq:j_matrix},
where the chemotactic contribution $C$ is non-symmetric. As a
consequence, the relationship between instability and
reactivity is no longer determined solely by the local
properties of $J^*$.

From Proposition \ref{prop:Jreactive_gen}, a necessary condition
for instability is the reactivity of the operator
$J^*-k^2C$, rather than the reactivity of $J^*$ itself.
Therefore, the destabilizing action of chemotaxis may compensate
for the lack of reactivity of the local kinetics. In particular,
whenever
$$|\beta|u^*>2\sqrt{D_uD_v},$$
the constraint requiring $J^*$ to be reactive disappears from
the necessary conditions for instability. In this regime,
chemotaxis may destabilize homogeneous equilibria that remain
non-reactive and asymptotically stable in the absence of
spatial transport.


This property marks a fundamental distinction between Turing and
chemotaxis-driven mechanisms of pattern formation. In
reaction--diffusion systems, instability necessarily originates
from a reactive equilibrium. In contrast, chemotaxis introduces
an additional route to instability through non-symmetric
transport, allowing pattern formation even when the local
kinetics lacks reactivity.

The geometric interpretation developed in Section
\ref{sec:chemo} and depicted in Figure  \ref{fig:primo_disegnino} further clarifies this mechanism. Indeed, when
$J^*$ is stable and non-reactive, the corresponding parameter
values belong to the half-plane $\eta<0$, where diffusion alone is
unable to induce instability. Nevertheless, sufficiently strong
chemotactic effects may move the system across the instability
boundary
$$\eta + g_u \theta =0,$$
thereby producing chemotaxis-driven instability despite the
absence of reactivity.

After introducing the discretization of the chemotactic term, the numerical examples will illustrate this phenomenon. While chemotaxis may simply modify existing Turing patterns, it may also generate spatial structures in parameter regimes where diffusion alone predicts neither instability nor pattern formation.

\section{Matrix-oriented discretization of the chemotactic term} 
\label{sec:matrix_chemotaxis}  
The theoretical results derived in the previous sections are independent of the particular spatial discretization. For the numerical experiments, however, we restrict our attention to rectangular domains $\Omega=[0,L_x]\times[0,L_y]$, discretized by a uniform Cartesian grid. This setting allows us to exploit the tensor-product structure of the discrete differential operator and to adopt a matrix-oriented formulation of the numerical method. 
The matrix-oriented approach for reaction--diffusion systems was introduced in \cite{DautiliaSguraSimoncini2020}, where the diffusion operator is written as a Kronecker sum and implicit time discretizations lead to Sylvester matrix equations involving one-dimensional differentiation matrices only. In that framework, the spatial discretization of the Laplacian is based on Extended Central Difference Formulas (ECDFs, \cite{amodio2005}), originally employed in the vector-based setting for the numerical approximation of Turing patterns in \cite{SGURA20124132}. These formulas provide high-order finite difference approximations of spatial derivatives while incorporating the homogeneous Neumann boundary conditions through suitable boundary closures. In the present work, for simplicity, we focus on the second-order ECDF approximation.
Our aim is to extend this matrix-oriented ECDF framework to reaction--diffusion systems with chemotaxis. This extension is not straightforward, since the chemotactic contribution is nonlinear, non-symmetric, and has the form of a cross-diffusion or taxis flux. Consequently, it cannot be incorporated directly into the standard matrix-oriented formulation developed for classical reaction--diffusion systems. To the best of our knowledge, matrix-oriented ECDF discretizations for reaction--diffusion systems with chemotactic cross-diffusion terms have not been systematically addressed. The goal of this section is therefore twofold. First, we recall the matrix-oriented form of the ECDF discretization of the Laplacian used in the implicit diffusion step. Second, we introduce a discretization of the chemotactic operator which is compatible with the same ECDF Laplacian. In particular, the proposed chemotactic discretization is constructed so that, in the constant-mobility case, it reduces exactly to the ECDF discrete Laplacian employed for diffusion.
Let $U(t), V(t) \in \mathbb{R}^{N_x\times N_y}$ denote the matrices collecting the nodal values of the two unknowns on the Cartesian grid of $\Omega=[0,L_x]\times[0,L_y]$. More precisely, $U_{ij}(t)\approx u(x_i,y_j,t), V_{ij}(t)\approx v(x_i,y_j,t)$, where $x_i=(i-1)h_x, \, i=1,\ldots,N_x,\, y_j=(j-1)h_y, \,j=1,\ldots,N_y.$ The one-dimensional ECDF matrices in the $x$- and $y$-directions are denoted by $T_x \in \mathbb{R}^{N_x\times N_x}, \, T_y\in\mathbb{R}^{N_y\times N_y}$. They include the homogeneous Neumann boundary closure. Hence, the two-dimensional ECDF Laplacian applied to the function $Z\in\mathbb{R}^{N_x\times N_y}$ is written in matrix form as 
\begin{equation} \label{eq:ecdf_laplacian_matrix} 
\Delta^{\mathrm{ECDF}} Z = T_x Z + Z T_y^T. 
\end{equation} 
Equivalently, using the column-wise vectorization operator, 
$$\operatorname{vec}(\Delta^{\mathrm{ECDF}} Z) = \left( I_{N_y}\otimes T_x + T_y \otimes I_{N_x} \right)\operatorname{vec}(Z),$$
one obtains the equivalent vector-based formulation. However, it is well known that the matrix-oriented formulation avoids the explicit construction of the full $N_x N_y\times N_x N_y$ matrix and preserves the two-dimensional structure of the grid.
The chemotactic term appearing in the first equation has the form $\nabla\cdot\bigl(u\nabla v\bigr)$. In order to obtain a discretization which is algebraically compatible with the ECDF Laplacian, we use the continuous identity \begin{equation} \label{eq:chemotaxis_product_identity} \nabla\cdot\bigl(u\nabla v\bigr) = u\Delta v + \nabla u\cdot\nabla v. \end{equation} 
This identity separates the Laplacian of $v$, which can be discretized by the same ECDF operator used for the diffusion, from the first-order gradient contributions. We therefore introduce first-derivative matrices 
$P_x\in\mathbb{R}^{N_x\times N_x}, P_y\in\mathbb{R}^{N_y\times N_y}$, constructed by using the same ghost-point closure as the ECDF second-derivative matrices. In the interior of the domain, the standard centered approximation is used,
$$\left(\frac{\partial z}{\partial x}\right)_i=\frac{z_{i+1}-z_{i-1}}{2h_x}, \qquad i=2,\ldots,N-1.$$ 
At the boundary nodes, the ECDF Neumann closure gives the corrected relations 
$$\left(\frac{\partial z}{\partial x}\right)_1 = \frac{2}{3h_x}(z_2-z_1), \qquad \left(\frac{\partial z}{\partial x}\right)_N = \frac{2}{3h_x}(z_N-z_{N-1}).$$ 
Thus, in one space dimension, the ECDF-compatible first-derivative operator is 
\begin{equation} 
\label{eq:D_ECDF} 
P_x = \frac{1}{2h_x} 
\begin{bmatrix} -\frac{4}{3} & \frac{4}{3} & 0 & \cdots & 0 \\ -1 & 0 & 1 & \ddots & \vdots \\ 0 & \ddots & \ddots & \ddots & 0 \\ \vdots & \ddots & -1 & 0 & 1 \\ 0 & \cdots & 0 & -\frac{4}{3} & \frac{4}{3} 
\end{bmatrix}
\end{equation} 
and analogously for $P_y$.
For a matrix grid function $Z$, the discrete first derivatives are then defined as 
\begin{equation} \label{eq:matrix_discrete_gradients} 
Z_x = P_x Z, \qquad Z_y = Z P_y^T. 
\end{equation} 
We define the ECDF-compatible discrete chemotactic operator by 
\begin{equation} 
\label{eq:chemotaxis_ECDF_standard}
\mathcal{C}^{\mathrm{ECDF}}(U,V) = U \odot \left(T_x V + V T_y^T\right) + \left[ (P_x U)\odot(P_x V) + (U P_y^T)\odot(V P_y^T) \right],
\end{equation} 
where $\odot$ denotes the Hadamard product, that is the elementwise multiplication of matrices. 
The main property of \eqref{eq:chemotaxis_ECDF_standard} is that it recovers exactly the ECDF Laplacian when the mobility is constant, that is, if we replace $U$ with $1$. In this sense, the chemotactic discretization is compatible with the discrete diffusion operator. 
\begin{proposition} 
\label{prop:ECDF_compatibility} 
Let $U = \mathbf{1}$ be the constant grid function in $\mathbb{R}^{N_x \times N_y}$. Then 
$$\mathcal{C}^{\mathrm{ECDF}}(\mathbf{1},V) = T_x V + V T_y^T = \Delta^{\mathrm{ECDF}} V.$$ 
\end{proposition} \begin{proof} 
Since the first-derivative matrices annihilate constant vectors, one has 
$$P_x \mathbf{1}_{N_x\times N_y} = P_x \left(\mathbf{1}_{N_x} \times \mathbf{1}_{N_y}  \right) = \left(P_x \mathbf{1}_{N_x}\right) \times \mathbf{1}_{N_y} = \mathbf{0} \times\mathbf{1}_{N_y} = \mathbf{0},$$
and analogously
$$\mathbf{1}_{N_x\times N_y}\, P_{y}^T= \left(\mathbf{1}_{N_x} \times \mathbf{1}_{N_y}  \right) \, P_{y}^T = \mathbf{1}_{N_x} \times \left(\mathbf{1}_{N_y} \, P_{y}^T \right) = \mathbf{1}_{N_y} \times \mathbf{0} = \mathbf{0}.$$
Therefore, for $U=\mathbf{1}_{N_x\times N_y}$, the two gradient terms in \eqref{eq:chemotaxis_ECDF_standard} vanish: $$(P_x U)\odot(P_x V) + (U P_y^T)\odot (V P_y^T) = 0. $$ Moreover, 
$$ U\odot\left(T_x V + V T_y^T\right) = T_x V + V T_y^T. $$ 
Hence 
$$ \mathcal{C}^{\mathrm{ECDF}}(\mathbf{1},V) = T_x V + V T_y^T = \Delta^{\mathrm{ECDF}} V,$$ 
which proves the claim. 
\end{proof} 
More generally, if $U=C_0$ is constant, then $$ \mathcal{C}^{\mathrm{ECDF}}(U,V) = C_0\left(T_x V + V T_y^T\right). $$ 
Thus, in the constant-mobility limit, the chemotactic operator reduces to the same ECDF Laplace operator used in the diffusion term. We now describe the corresponding time discretization. We use an IMEX Euler scheme in which diffusion is treated implicitly, while reaction and chemotaxis are treated explicitly. Given $U^n,V^n$ at time $t_n=n h_t$, define 
\begin{equation} 
\label{eq:rhs_u_matrix_chem} 
R_u^n = U^n + h_t\,F(U^n,V^n) - h_t\,\beta\,\mathcal{C}^{\mathrm{ECDF}}(U^n,V^n), 
\end{equation} 
and 
\begin{equation}
\label{eq:rhs_v_matrix_chem} 
R_v^n = V^n + h_t\,G(U^n,V^n), \end{equation} 
where $F$ and $G$ denote the componentwise evaluations of the reaction terms. The new iterates $U^{n+1}$ and $V^{n+1}$ are obtained by solving the Sylvester equations \begin{equation} \label{eq:sylvester_u_chem} \left(I_{N_x}-h_t\,D_u T_x\right)U^{n+1} + U^{n+1} \left(-h_t\,D_u T_y^T\right) = R_u^n, \end{equation} 
and 
\begin{equation} \label{eq:sylvester_v_chem} \left(I_{N_x}-h_t\,D_v T_x\right)V^{n+1} + V^{n+1} \left(-h_t\,D_v T_y^T\right) = R_v^n. \end{equation} 
Therefore, the inclusion of chemotaxis does not destroy the matrix-oriented structure of the implicit diffusion step. The chemotactic term only modifies the explicit right-hand side \eqref{eq:rhs_u_matrix_chem}, while the implicit part is solved by the same Sylvester strategy used for pure reaction--diffusion systems. Since the coefficient matrices in \eqref{eq:sylvester_u_chem}--\eqref{eq:sylvester_v_chem} are constant in time, their spectral decompositions can be precomputed. For example, setting 
$$ A_u = I_{N_x}-h_t\,D_u T_x, \qquad B_u = -h_t\,D_u T_y^T, $$ 
and assuming 
$$ A_u = Q_{xu}\,\Lambda_{xu}\,Q_{xu}^{-1}, \qquad B_u = Q_{yu}\,\Lambda_{yu}\,Q_{yu}^{-1}, $$ 
the solution of \eqref{eq:sylvester_u_chem} is computed by 
$$ \widehat{Q}_u^n = Q_{xu}^{-1}Q_u^nQ_{yu}, $$ 
followed by the componentwise solve $$ \widehat{U}^{n+1}_{ij} = \frac{\widehat{Q}^{\,n}_{u,ij}} {(\Lambda_{xu})_{ii}+(\Lambda_{yu})_{jj}}, $$ 
and the projection back to the physical basis, 
$$ U^{n+1} = Q_{xu}\widehat{U}^{n+1}Q_{yu}^{-1}. $$ 
The same procedure is applied to the equation for $V^{n+1}$. This reduced implementation preserves the computational advantages of the matrix-oriented method while including the chemotactic contribution through the explicit term $\mathcal{C}^{\mathrm{ECDF}}(U^n,V^n)$. 

\begin{remark} 
The discretization \eqref{eq:chemotaxis_ECDF_standard} is designed to be algebraically compatible with the ECDF Laplacian. This should be distinguished from flux-based finite-volume or upwind discretizations of chemotaxis, which are often introduced to enforce discrete mass conservation, positivity, or robustness in taxis-dominated regimes. Such conservative approaches are widely used for Keller--Segel and more general taxis--diffusion--reaction systems; see, for instance, \cite{GerischChaplain2006,Filbet2006,Saito2007,EpshteynKurganov2008, ChertockEpshteynHuKurganov2018}. In contrast, the present construction prioritizes compatibility with the ECDF matrix Laplacian used in the implicit diffusion solver, as expressed by the identity $\mathcal{C}^{\mathrm{ECDF}}(\mathbf{1},V)=\Delta^{\mathrm{ECDF}} V. $ 
\end{remark}

\section{Numerical examples}
\label{sec:numerical_examples}
In this section, we present a series of numerical experiments aimed at 
illustrating how chemotaxis can enhance, or modify the geometry of
classical diffusion-driven Turing patterns.  
We first focus on reaction-diffusion models for which  the Jacobian $J^*$
evaluated at equilibrium is reactive, i.e., the homogeneous steady state is stable and \emph{reactive}, exhibiting short-time amplification of perturbations. In these conditions, patterns emerge through the standard Turing mechanism driven solely by diffusion.  
Our goal is to investigate whether, and to what extent, these diffusion-induced 
patterns persist when a chemotactic term is added, and how the interplay 
between diffusion, reaction kinetics, and chemotaxis modifies the geometry and intensity of the resulting 
spatial structures.

\subsection{Chemotaxis-enhanced Turing patterns:
the Schnakenberg model}

To illustrate our analysis, we consider the classical Schnakenberg reaction-diffusion system \cite{Schnakenberg1979,Murray_book2} 
\begin{equation}
\label{schnak_model}
\begin{aligned}
  \dfrac{\partial u}  {\partial t} & = D_u\, \Delta u + \gamma \left(a - u + u^2 v\right),\\
\dfrac{\partial v}{\partial t} & =  D_v \,\Delta v + \gamma \left(b - u^2 v\right),
\end{aligned}
\end{equation}
a prototypical two-species model widely employed to study pattern formation in chemical and biological contexts. The system is posed on a bounded domain $\Omega \subset \mathbb{R}^2$ and complemented with homogeneous Neumann boundary conditions.

However, in the present work we are interested in understanding how the geometry and stability of Turing patterns are modified when a chemotactic mechanism is introduced. To this end, we extend system \eqref{schnak_model} by adding a standard Keller-Segel type flux, leading to
\begin{equation}
    \label{schnak_chemo}
    \begin{aligned}
  \dfrac{\partial u}  {\partial t} & = D_u\, \Delta u -\beta \, \nabla \cdot(u\nabla v)+ \gamma \left(a - u + u^2 v\right),\\
\dfrac{\partial v}{\partial t} & =  D_v \,\Delta v + \gamma \left(b - u^2 v\right),
\end{aligned}
\end{equation}
where the parameter $\beta>0$ denotes the chemotactic sensitivity of species $u$ 
towards gradients of $v$.
For any choice of positive parameters $a$ and $b$, the model admits a unique spatially uniform steady state
\begin{equation}
\label{eq:schank_pe}
P_e = (u^*,v^*) = \left(a+b,\frac{b}{(a+b)^2} \right),
\end{equation}
whose stability depends on the interplay between the parameters 
$a$, $b$.

We first observe that the kinetics satisfies 
\begin{equation}
\label{eq:gu_schnak}
g_u(u^*,v^*) = -\frac{2b}{a+b}<0
\end{equation}
at the homogeneous equilibrium $P_e$ \eqref{eq:schank_pe}.  
This corresponds to the configuration represented in the right panel of 
Figure \ref{fig:primo_disegnino}. As already discussed in Section \ref{sec:chemo}, increasing the parameter $\beta$ shifts 
the system toward the region of 
linear stability.  
Therefore, sufficiently large values of $\beta$ can push the system into a 
stable regime, suppressing the diffusion-driven instability.

The following theorem identifies the exact critical value of the chemotactic 
sensitivity $\beta^*$ at which this stabilizing effect cancels the Turing 
instability.

\begin{theorem}
If 
\begin{equation}
    \label{schnak_betac}
    \beta < \beta^* = -D_u \frac{(a+b)^2}{2b} + D_v \frac{b-a}{2b(a+b)} - \sqrt{D_u D_v} \frac{a+b}{b}
\end{equation}
then the spatially homogeneous equilibrium $P_e$ \eqref{eq:schank_pe} of model \eqref{schnak_chemo} undergoes diffusion-chemotaxis-driven instability.
\end{theorem}

The numerical experiments presented below show that chemotaxis affects both 
the intensity and spatial distribution of the spots, with differences whose 
magnitude depends on the value of the chemotactic sensitivity.

\subsubsection{Spot-like patterns} 
\begin{figure}[ht] 
\centering \includegraphics[width=0.55\linewidth]{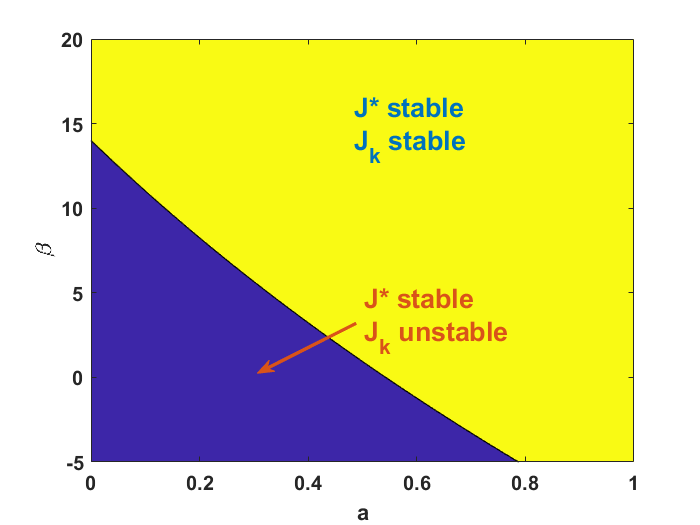} 
\caption{ Bifurcation diagram in the $(a,\beta)$-parameter plane for the chemotaxis-extended Schnakenberg model with $D_v = 100$. The yellow region corresponds to parameter values for which the homogeneous equilibrium is linearly stable with respect to spatial perturbations, while the blue region denotes diffusion--chemotaxis-driven instability. Since $g_u(P_e)<0$, positive values of $\beta$ have a stabilizing effect, whereas negative values of $\beta$ enlarge the instability region. The arrow indicates the direction in which the system moves when $\beta$ is decreased.} \label{fig:schnak_bifurcation} 
\end{figure}

We first consider a parameter regime for which the classical Schnakenberg reaction--diffusion system exhibits a diffusion-driven Turing instability leading to spot-like patterns \cite{AMS2023,BMS2021}. The kinetic and diffusive parameters are chosen as in the VisualPDE experiments of Walker et al.~\cite{walker2023visualpde}, namely, 
$$D_u=1, \qquad D_v=100, \qquad a=0.01, \qquad b=2 .$$ 
The scaling parameter is fixed at $\gamma=1000$. The system is solved on the square domain $\Omega=[0,2]\times[0,2],$ discretized by a uniform Cartesian grid with $N_x=N_y=100$ nodes on each spatial direction. The time integration is performed by the IMEX Euler scheme described in Section~\ref{sec:matrix_chemotaxis}, with time step $h_t=10^{-4}$, up to the final time $T_f=50$. Homogeneous Neumann boundary conditions are imposed on both variables. All simulations are initialized by perturbing the spatially homogeneous equilibrium with the same small-amplitude random perturbation. As observed in \eqref{eq:gu_schnak}, $g_u(P_e)=<0.$ Thus, according to the geometric interpretation discussed in Section~\ref{sec:chemo}, positive values of the chemotactic sensitivity $\beta$ move the system towards the stable region, whereas negative values of $\beta$ reinforce the instability. This behavior is summarized in the bifurcation diagram shown in Figure~\ref{fig:schnak_bifurcation}. 

\begin{figure}[t] 
\centering 
\includegraphics[width=0.32\linewidth]{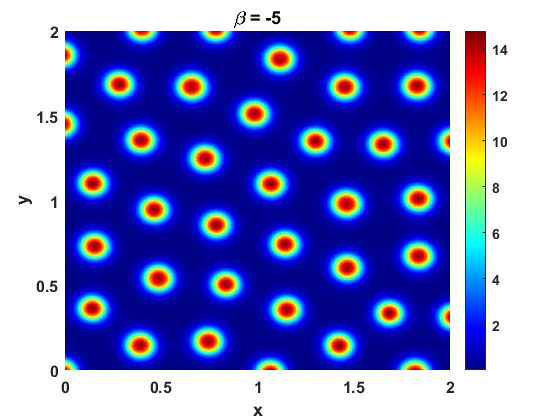} 
\includegraphics[width=0.32\linewidth]{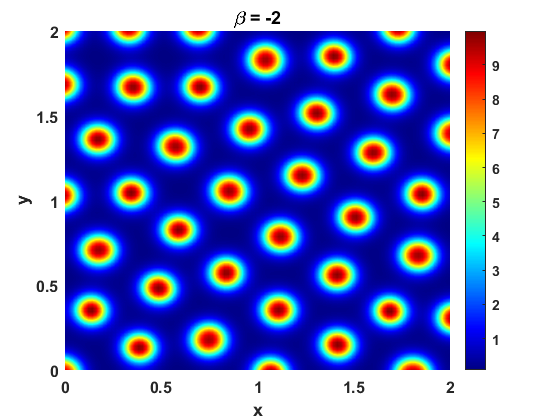} 
\includegraphics[width=0.32\linewidth]{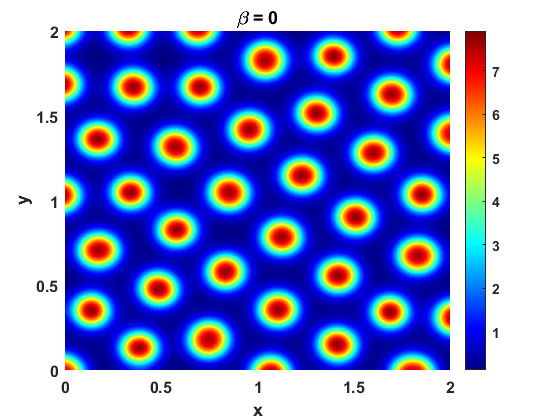} 
\includegraphics[width=0.32\linewidth]{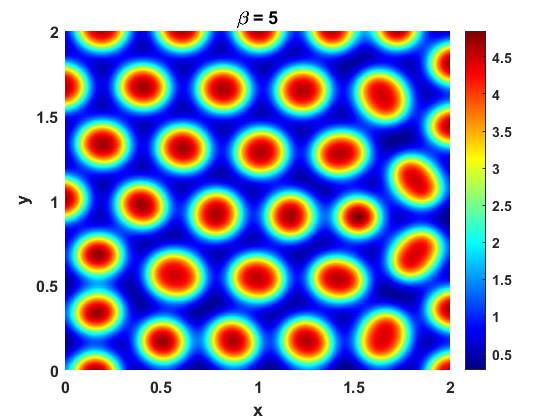} 
\includegraphics[width=0.32\linewidth]{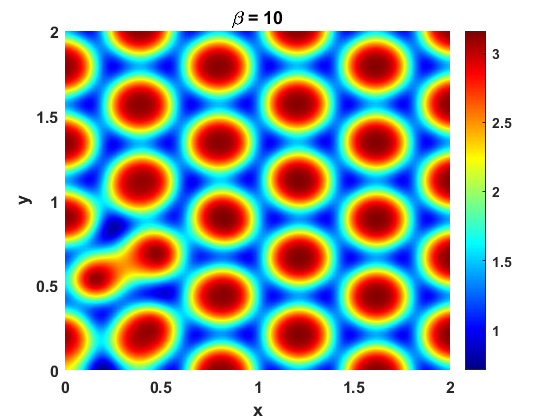} 
\includegraphics[width=0.32\linewidth]{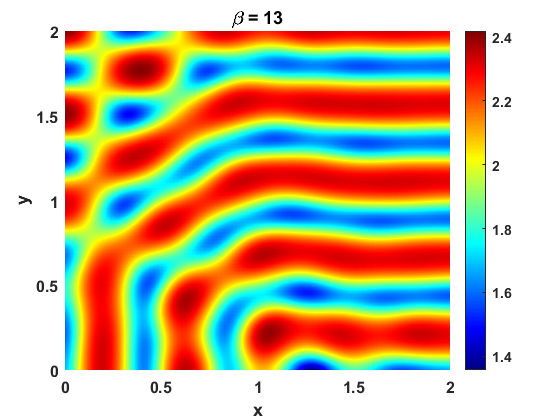} 
\caption{Effect of the chemotactic sensitivity $\beta$ on the spot-like Turing pattern of the Schnakenberg model. The panels show $u(x,y,T_f)$ at $T_f=50$ for different values of $\beta$. For $\beta=0$, the model reduces to the classical reaction--diffusion Schnakenberg system and produces a regular spot-like Turing pattern. Negative values of $\beta$ enhance the instability, yielding sharper spots with larger peak values. Positive values of $\beta$ have a stabilizing effect: the spots become wider, their amplitude decreases, and the solution progressively evolves towards mixed spot--stripe configurations before approaching the homogeneous steady state near the critical threshold $\beta^*$. Each panel is plotted with its own color scale in order to highlight the spatial morphology; the colorbars also show the variation in the amplitude of $u$.} \label{fig:schnak_spots} 
\end{figure}

The yellow region corresponds to parameter values for which the homogeneous equilibrium is stable both in the absence and in the presence of spatial perturbations, namely $J^*$ is stable and $J_k$ is stable for all admissible wavenumbers. The blue region corresponds instead to diffusion--chemotaxis-driven instability: $J^*$ remains stable, but $J_k$ becomes unstable for some $k \neq 0$ Since $g_u(P_e)<0$, increasing $\beta$ shifts the system towards the stable region, while decreasing $\beta$  moves it deeper into the instability region. The corresponding spatial patterns are displayed in Figure~\ref{fig:schnak_spots}. 

For $\beta = 0$, system \eqref{schnak_chemo} reduces to the classical reaction--diffusion Schnakenberg model. In this case, diffusion destabilizes the homogeneous equilibrium and the solution develops a stationary spot-like Turing pattern. The pattern consists of localized regions of high $u$-concentration arranged in an approximately regular configuration. When $\beta<0$ , the chemotactic contribution enhances the instability. The spots become more localized and their peak values increase, as indicated by the larger range of the corresponding color scale. This confirms that, in the regime $g_u<0$, negative chemotaxis reinforces the aggregation mechanism responsible for the formation of localized high-$u$ structures. In contrast, for $\beta>0$, chemotaxis progressively counteracts the diffusion-driven instability. As $\beta$  increases, the spots become wider and less pronounced, and their maximum values decrease. For intermediate positive values of $\beta$, the regular spot lattice is lost and the solution develops mixed spot--stripe configurations. For values of $\beta$  close to, or larger than, the critical threshold $\beta^* \approx $ in \eqref{schnak_betac}, the spatial modulation is strongly damped and the solution approaches the homogeneous steady state. Therefore, in this example chemotaxis does not generate a new instability mechanism, since the system is already Turing unstable for $\beta=0$. Rather, it modulates the pre-existing diffusion-driven pattern: negative chemotaxis strengthens the instability, while positive chemotaxis weakens it and may eventually suppress pattern formation.

\subsubsection{Stripe-like patterns}
We now consider a second parameter regime of the Schnakenberg model, again inspired by the VisualPDE experiments, in which the classical reaction--diffusion system produces stripe-like patterns. All parameters are kept fixed as in the previous experiment, except for the diffusion coefficient of the second species, which is now set to $D_v=30$. The remaining numerical parameters are unchanged.  

\begin{figure}[t] \centering \includegraphics[width=0.55\linewidth]{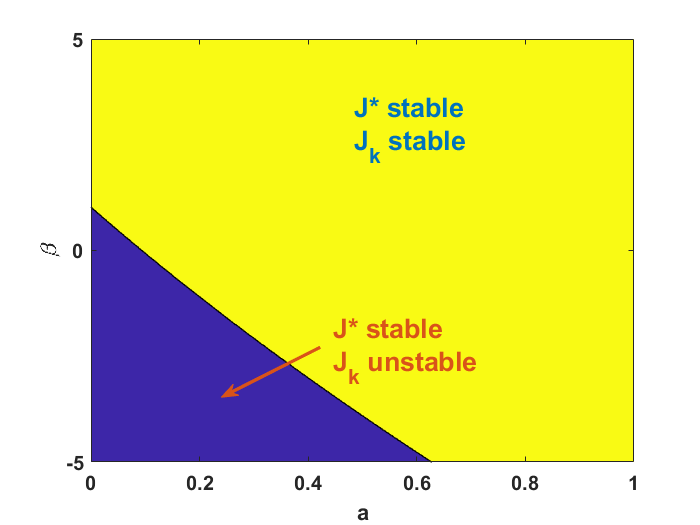}
\caption{Bifurcation diagram in the $(a,\beta)$-parameter plane for the chemotaxis-extended Schnakenberg model with $D_v=30$. The yellow region denotes parameter values for which the homogeneous equilibrium is stable with respect to spatial perturbations, while the blue region corresponds to diffusion--chemotaxis-driven instability. Since $g_u(P_e)<0$, negative values of $\beta$ enhance the instability, whereas positive values of $\beta$ move the system towards the stable region.} \label{fig:schnak_bifurcation_stripes} \end{figure}

The corresponding bifurcation diagram in the $(a,\beta)$-parameter plane is shown in Figure~\ref{fig:schnak_bifurcation_stripes}. 

\begin{figure}[t] \centering 
\includegraphics[width=0.32\linewidth]{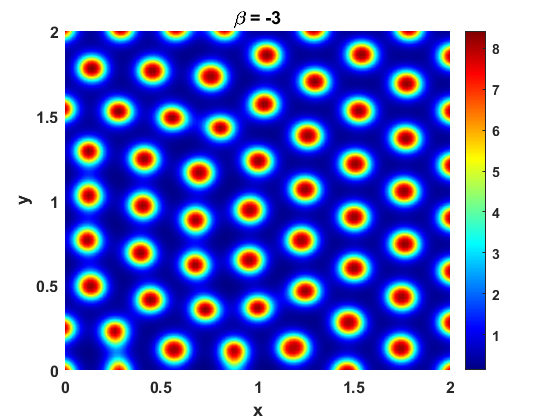} 
\includegraphics[width=0.32\linewidth]{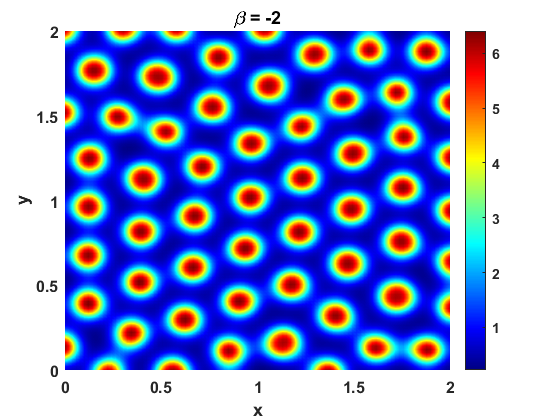} 
\includegraphics[width=0.32\linewidth]{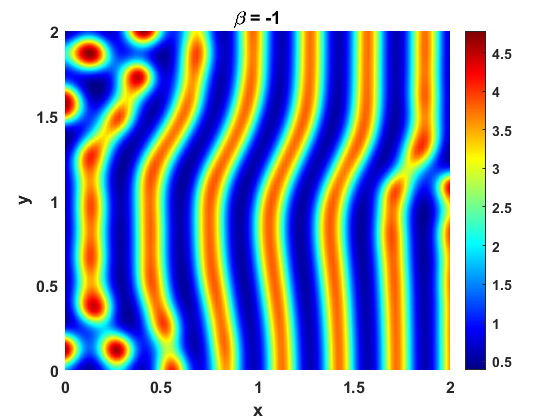} 
\includegraphics[width=0.32\linewidth]{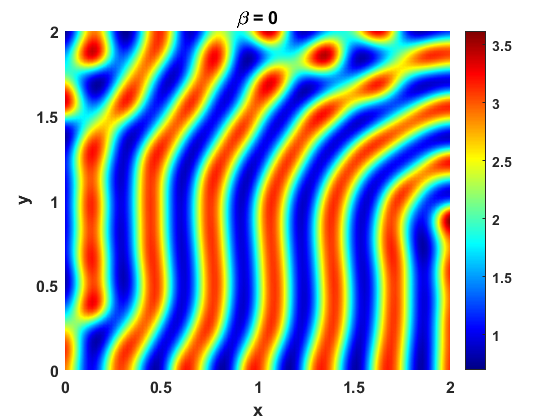} 
\includegraphics[width=0.32\linewidth]{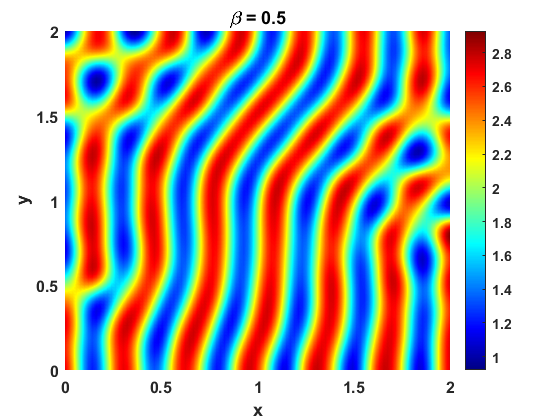} 
\includegraphics[width=0.32\linewidth]{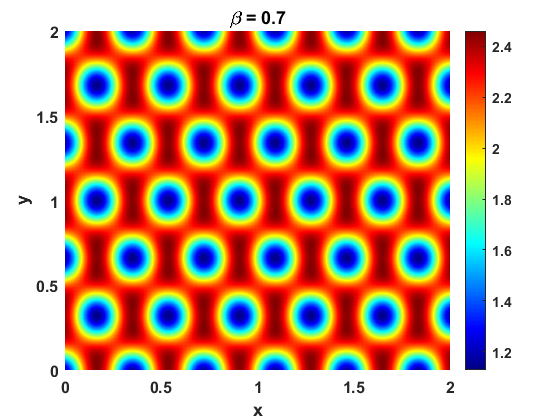} 
\caption{Effect of the chemotactic sensitivity $\beta$ on stripe-like Turing patterns of the Schnakenberg model with $D_v=30$. The panels show $u(x,y,T_f)$ at $T_f=50$. For $\beta=0$, the pure reaction--diffusion system produces a stripe-like pattern. Negative values of $\beta$ enhance the instability and lead to high-amplitude localized spots. Positive values of $\beta$ reduce the pattern amplitude but also induce a transition from stripes to lower-amplitude spot-like or cellular structures before the spatial modulation is eventually suppressed. Each panel is plotted with its own color scale.} \label{fig:schnak_stripes} 
\end{figure}

As in the previous case, since $g_u(P_e)<0$, the role of the chemotactic parameter is the same as in the spot-like experiment: negative values of $\beta$ move the system deeper into the instability region, whereas positive values of $\beta$ shift the system towards the stability boundary. The associated spatial patterns are displayed in Figure~\ref{fig:schnak_stripes}. 

For $\beta=0$, the system reduces to the classical reaction--diffusion Schnakenberg model and develops a stripe-like Turing pattern, consisting of alternating bands of high and low $u$-concentration. For $\beta<0$, chemotaxis enhances the instability. The amplitude of the solution increases and the stripe structure progressively breaks into localized high-amplitude spots. In contrast, for $\beta>0$, the overall amplitude decreases, consistently with the stabilizing effect predicted by the bifurcation diagram. Nevertheless, before reaching the homogeneous steady state, the stripe pattern undergoes a morphological transition towards lower-amplitude spot-like or cellular structures. This experiment shows that, in the stripe-forming regime, chemotaxis acts not only as an amplitude-modulating mechanism but also as a morphology-selecting mechanism. Spot-like structures may appear on both sides of $\beta=0$, but with different dynamical meaning: for $\beta<0$ they result from a reinforcement of the instability, whereas for $\beta>0$ they arise during the weakening of the stripe instability and precede the suppression of the pattern near the critical chemotactic threshold.

\subsection{Chemotaxis-generated patterns:
a predator--prey model}

The Schnakenberg example illustrates how chemotaxis modifies pattern-forming mechanisms that are already present in a purely reaction--diffusion framework. In that case, spatial structures arise through the classical Turing instability and chemotaxis acts either as a stabilizing or a destabilizing factor depending on the sign of the coupling.

We now turn to a fundamentally different situation, which is more directly related to the main message of this work. The theoretical results derived in Section~\ref{sec:chemo} show that reactivity of the kinetic Jacobian is not a necessary condition for the onset of chemotaxis-driven instability. To illustrate this phenomenon, we consider a predator--prey system for which diffusion alone cannot generate patterns and whose kinetic Jacobian is non-reactive. In this setting, any observed spatial structure must originate from the chemotactic mechanism itself rather than from an underlying Turing instability. This scenario is not limited to predator--prey systems. Similar non-reactive kinetics arise, for instance, in the MOMOS model for soil organic carbon dynamics studied in \cite{Monti_MATCOM2025}, motivating the present investigation into chemotaxis-driven instabilities beyond the reactive setting.

To this end, we consider a predator--prey model inspired by the prey--taxis system recently proposed by Wang et al.~\cite{Wang2026}, where directed predator movement is incorporated into a reaction--diffusion framework through a prey-taxis mechanism.

The model reads
\begin{equation}
\label{eq:ppHT}
\begin{aligned}
\frac{\partial u}{\partial t}
&=
D_u\Delta u
-\beta\nabla\!\cdot(u\nabla v)
+\gamma\left(
-d\,u
+\frac{cuv}{(1+av)(1+bu)}
\right),
\\[1ex]
\frac{\partial v}{\partial t}
&=
D_v\Delta v
+\gamma\left(
v(1-v)
-\frac{uv}{(1+av)(1+bu)}
\right),
\end{aligned}
\qquad x\in\Omega,\; t>0,
\end{equation}
where $\Omega \subset \mathbb{R}^2$ and supplemented with homogeneous Neumann boundary conditions.

Here $u=u(x,y,t)$ and $v=v(x,y,t)$ denote the predator and prey densities, respectively, while $D_u,D_v>0$ are the corresponding diffusion coefficients.

The term $-\beta \nabla\cdot(u\nabla v)$ describes the directed movement of predators in response to the spatial distribution of prey. As in classical chemotaxis models, predators may actively adjust their movement according to resource availability. When $\beta>0$, predators move preferentially toward regions with higher prey density, modeling a foraging strategy driven by food attraction. Conversely, when $\beta<0$, predators move away from regions of increasing prey density. Such a retreat behavior may arise as a consequence of prey defense mechanisms or other ecological factors that make prey-rich regions less favorable for predators. Therefore, the parameter $\beta$ determines both the intensity and the direction of predator movement induced by prey density gradients.

The reaction kinetics combine logistic prey growth with a Crowley--Martin functional response \cite{Crowley1989}. The prey population grows according to the logistic law
$v(1-v)$, while predation is described by
$\displaystyle \frac{uv}{(1+av)(1+bu)}$,
which accounts simultaneously for prey handling effects and predator interference. The parameter $a>0$ is associated with prey handling, whereas $b>0$ measures the strength of predator interference.

The predator population experiences a constant mortality rate $d>0$ and increases through prey consumption with conversion efficiency $c>0$. Finally, $\gamma >0$ is the classical scaling parameter introduced through nondimensionalization. Following the interpretation of Murray \cite{Murray_book2}, increasing $\gamma$ corresponds to increasing the domain size and, equivalently, to strengthening the reaction kinetics relative to diffusion.

The interplay between nonlinear predator--prey interactions, diffusion, and directed movement provides a natural mechanism for spatial self-organization. In particular, the prey-taxis coefficient $\beta$ acts as the bifurcation parameter regulating the transition from homogeneous equilibria to heterogeneous spatial patterns.

For the parameter values considered throughout this work, the system admits a unique positive spatially homogeneous equilibrium $P_e=(u^*,v^*),$
where $(u^*,v^*)$ denotes the unique positive solution of the nonlinear algebraic system obtained by setting the reaction terms equal to zero.

We first observe that
$$g_u(u^*,v^*)
=
-\gamma
\frac{v^*}{(1+av^*)(1+bu^*)^2}<0.$$

This corresponds to the configuration represented in the right panel of Figure~\ref{fig:primo_disegnino}. 


The following result identifies the critical value of the chemotactic sensitivity above which the homogeneous equilibrium remains spatially stable, and below which diffusion--chemotaxis driven instability occurs. \begin{theorem} 
Let $P_e=(u^*,v^*)$ be a positive spatially homogeneous equilibrium of model~\eqref{eq:ppHT}, and assume that $J^*$ is asymptotically stable, namely $$ \operatorname{tr}(J^*)<0, \qquad \det(J^*)=f_u g_v-f_v g_u>0. $$  
Then, the equilibrium $P_e$ undergoes diffusion--chemotaxis-driven instability whenever \begin{equation} \label{eq:pp_beta_crit} \beta < \beta^* := \frac{ 2\sqrt{D_uD_v\left(f_u g_v-f_v g_u\right)} -\left(D_u g_v+D_v f_u\right) } {g_u u^*}, \end{equation} where all the partial derivatives are evaluated at $P_e$. 
\end{theorem}
\begin{proof} 
The instability condition recalled in Section~\ref{sec:chem_insta} requires $$ D_u g_v+D_v f_u+\beta g_u u^* > 2\sqrt{D_uD_v\left(f_u g_v-f_v g_u\right)}. $$ Since $g_u(P_e)<0$ and $u^*>0$, the coefficient $g_u u^*$ is negative. Solving the previous inequality with respect to $\beta$ gives $$ \beta < \frac{ 2\sqrt{D_uD_v\left(f_u g_v-f_v g_u\right)} -\left(D_u g_v+D_v f_u\right) } {g_u u^*} = \beta^*. $$ This proves the claim. 
\end{proof}

If the homogeneous equilibrium lies strictly outside the purely diffusive Turing regime, then 
$$D_u g_v+D_v f_u < 2\sqrt{D_uD_v\left(f_u g_v-f_v g_u\right)}.$$ Therefore, the numerator in \eqref{eq:pp_beta_crit} is positive. On the other hand, since $g_u(P_e)<0$ and $u^*>0$, the denominator $g_u u^*$ is negative. Hence, $$ \beta^*<0. $$ As a consequence, diffusion--chemotaxis-driven instability can occur only for sufficiently negative values of the prey-taxis coefficient, namely for $\beta<\beta^*$. From the biological viewpoint, this corresponds to a repulsive prey-taxis mechanism, whereby predators move away from regions of high prey density.

\begin{figure}[t] \centering 
\includegraphics[width=0.45\linewidth]{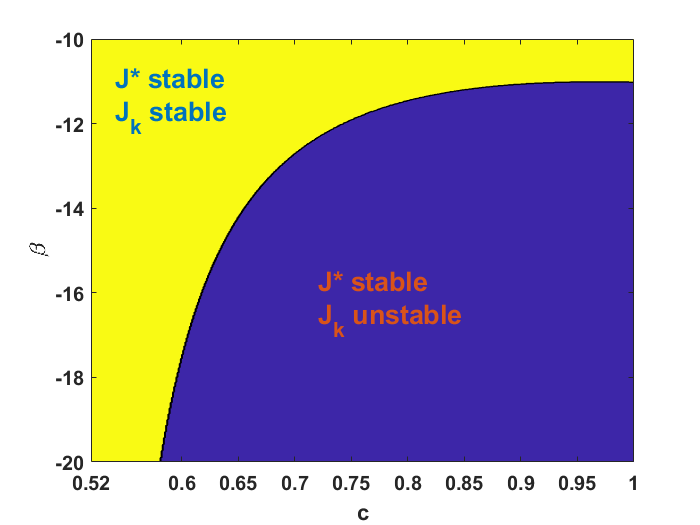}
\includegraphics[width=0.45\linewidth]{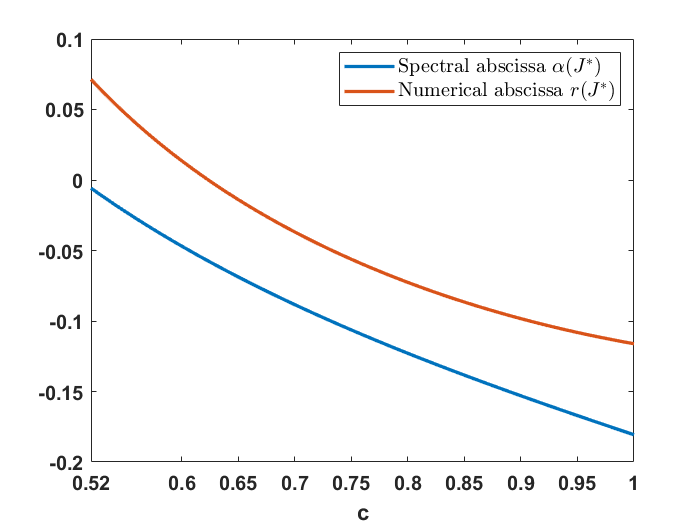}
\caption{Bifurcation and reactivity analysis for the predator--prey model~\eqref{eq:ppHT}. 
Left: bifurcation diagram in the $(c,\beta)$-parameter plane. The curve represents the critical threshold $\beta=\beta^*(c)$ defined in~\eqref{eq:pp_beta_crit}; the blue region below the curve corresponds to diffusion--chemotaxis-driven instability, while the yellow region above the curve corresponds to spatial stability of the homogeneous coexistence equilibrium. The interval starts at $c=0.52$, slightly above the coexistence threshold $c^*=d(1+a)=0.51$. 
Right: spectral abscissa $\alpha(J^*)=\max_i\Re(\lambda_i(J^*))$ and numerical abscissa $r(J^*)=\lambda_{\max}(H(J^*))$, as functions of $c$. The spectral abscissa remains negative throughout the considered interval, showing that the coexistence equilibrium is asymptotically stable. The numerical abscissa is positive only in a small interval close to the lower coexistence threshold, while it becomes negative for larger values of $c$. In particular, the value $c=1$, used in the pattern simulations, lies in the stable and non-reactive regime.
} \label{fig:pp_bifurcation}
\end{figure}

The numerical experiments presented below illustrate this mechanism through a bifurcation diagram in the $(c,\beta)$-parameter plane and the corresponding stationary spatial patterns, confirming that chemotaxis alone can trigger pattern formation even in the absence of reactivity and, thus, classical Turing instability. For the numerical simulations, we use the parameter values considered in~\cite{Wang2026}, 
\begin{equation}
\label{eq:pp_param}
D_u=0.55,\qquad D_v=2,\qquad a=0.7,\qquad b=1,\qquad d=0.3,
\end{equation}
and we set $\gamma=500.$ 
For these parameters, the positive coexistence equilibrium exists for
$c>d(1+a)=0.51.$
Therefore, in the diagrams reported in Figure \ref{fig:pp_bifurcation}  we consider the interval $c\in[0.52,1]$, which lies inside the coexistence regime. 
All simulations are performed on the square domain $\Omega=[0,2]\times[0,2],$
discretized by a uniform Cartesian grid with $N_x=N_y=100$ nodes in each spatial direction. Homogeneous Neumann boundary conditions are imposed on both variables. Time integration is carried out by the IMEX Euler scheme described in Section~\ref{sec:matrix_chemotaxis}, with time step
$h_t=10^{-4}$. The initial condition is obtained by perturbing the spatially homogeneous coexistence equilibrium with a small-amplitude random perturbation. Unless otherwise stated, the solution is evolved up to the final time $T_f=1000$.

Figure \ref{fig:pp_bifurcation} summarizes both the diffusion--chemotaxis instability threshold and the local stability/reactivity properties of the coexistence equilibrium. The left panel shows the bifurcation diagram in the $(c,\beta)$-plane, while the right panel displays the spectral abscissa $\alpha(J^*)$ and the numerical abscissa $r(J^*)$ of the kinetic Jacobian as functions of $c$. 
The left panel shows that the line $\beta=0$ lies in the stable region for the whole interval of $c$ considered. Hence, diffusion alone does not generate a Turing instability in this parameter range. Spatial instability appears only for sufficiently negative values of the chemotactic sensitivity, namely below the threshold $\beta=\beta^*(c)$.
The right panel confirms that the coexistence equilibrium is asymptotically stable throughout the considered interval, since $(\alpha(J^*)<0$. The numerical abscissa $r(J^*)$, however, is positive only in a small range close to the lower coexistence threshold and becomes negative for larger values of $c$. In particular, the value $c=1$, used in the pattern simulations below, lies well inside the non-reactive regime. This choice is therefore especially relevant for the purpose of this work: it allows us to test whether chemotaxis can generate spatial patterns when the homogeneous equilibrium is stable and non-reactive.

For $c=1$, the positive homogeneous equilibrium is
$P_e=(u^*,v^*) \simeq (0.6327,0.7454).$
The kinetic Jacobian evaluated at $P_e$, up to the common scaling factor $\gamma$, is
$$J^* \simeq
\begin{pmatrix}
-0.1163 & 0.1673\\
-0.1837 & -0.6580
\end{pmatrix}.
$$
Its eigenvalues are $\lambda_1\simeq -0.1807,
\lambda_2\simeq -0.5936,$
so that $\alpha(J^*)<0$. Moreover, the eigenvalues of the Hermitian part
$H(J^*)=\frac{J^*+(J^*)^T}{2}$
are
$\mu_1\simeq -0.6582, \mu_2\simeq -0.1161.$
Therefore, $r(J^*)=\lambda_{\max}(H(J^*))\simeq -0.1161<0,$
showing that the homogeneous coexistence equilibrium used in the simulations is non-reactive. The scaling factor $\gamma=500$ multiplies these eigenvalues by the same positive factor and hence does not change their sign.

\begin{figure}[t] 
\centering 
\includegraphics[width=0.32\linewidth]{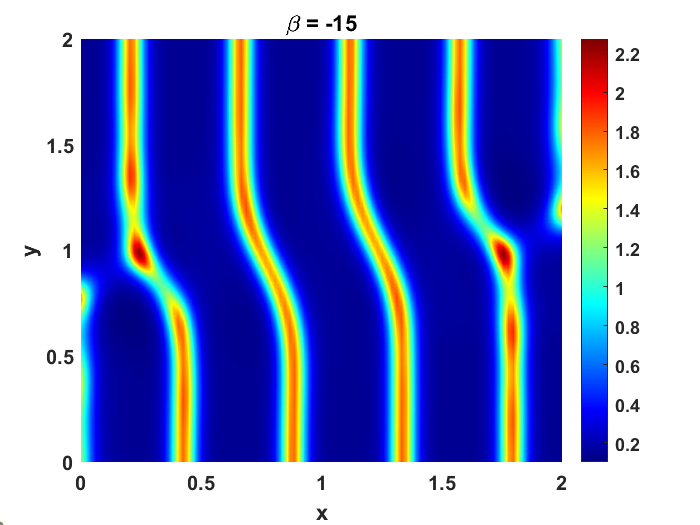}
\includegraphics[width=0.32\linewidth]{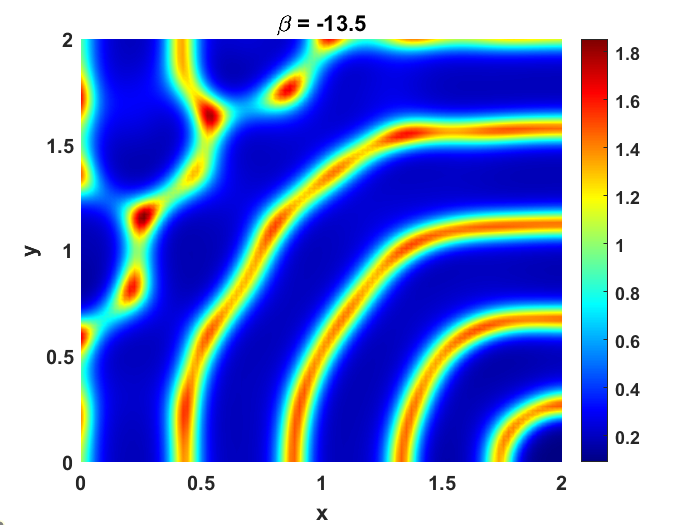}
\includegraphics[width=0.32\linewidth]{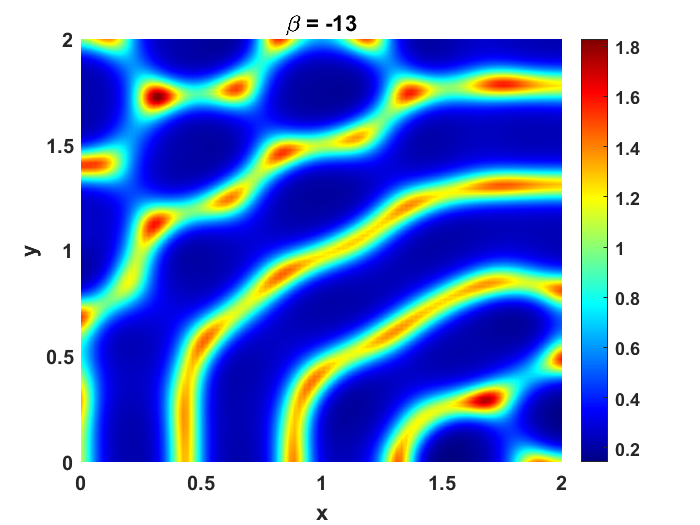}
\includegraphics[width=0.32\linewidth]{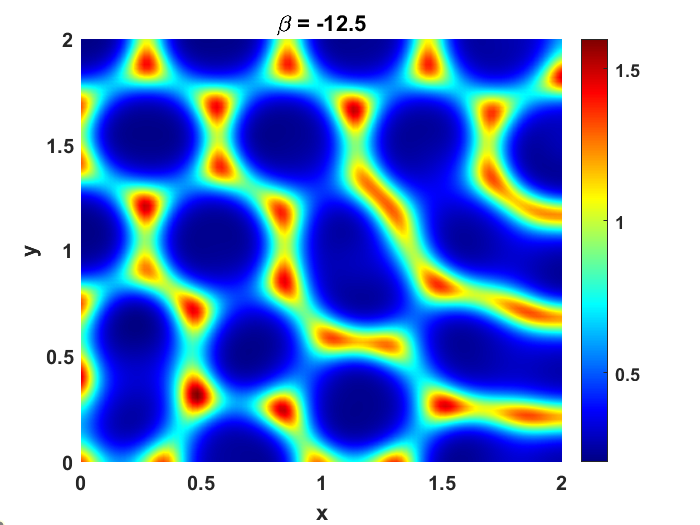}
\includegraphics[width=0.32\linewidth]{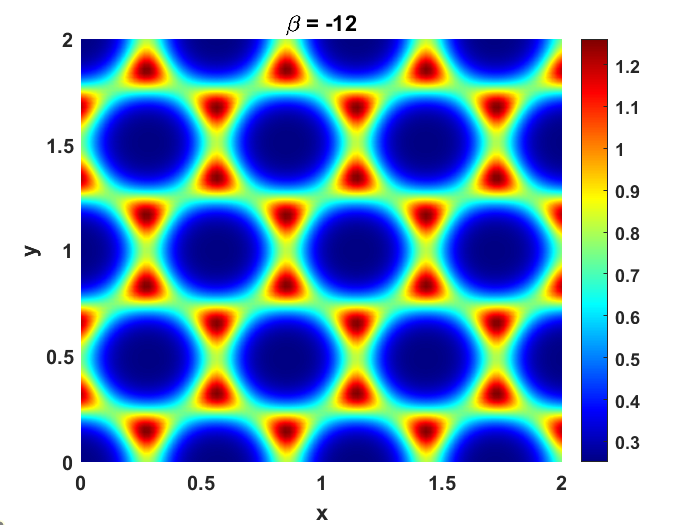}
\includegraphics[width=0.32\linewidth]{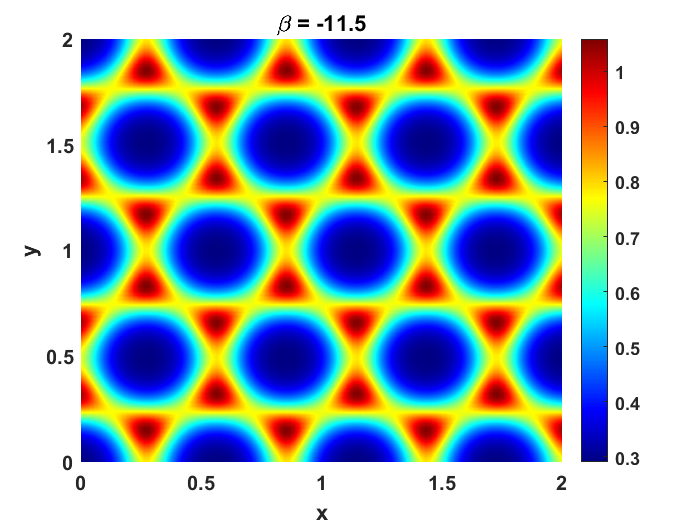}
\caption{
Effect of the chemotactic sensitivity $\beta$ on the stationary patterns of the predator--prey model \eqref{eq:ppHT}.
The panels show the predator density $u(x,y,T_f)$ at $T_f=1000$ for different values of $\beta$, with $c=1$ and the remaining parameters fixed as in \eqref{eq:pp_param}.
As $\beta$ approaches the stability threshold $\beta^*$ from below, the pattern amplitude decreases and the solution organizes into nearly regular hexagonals. For more negative values of $\beta$, corresponding to parameter regimes located further inside the diffusion--chemotaxis instability region, the amplitude increases and the hexagonal-like structures progressively merge into stripe-like patterns.
The figure therefore reveals both an amplitude modulation and a morphology transition from hexagons to stripes as the system moves away from the stability boundary.
Each panel is plotted with its own color scale in order to highlight the spatial morphology; the colorbars also show the variation in the amplitude of $u$.} \label{fig:pp_WLL_pattern} 
\end{figure}

For the same parameter values, the critical chemotactic threshold defined in \eqref{eq:pp_beta_crit} is
$\beta^* \approx -11.02.$
Thus, the values $\beta=-15,\,-13.5,\,-13,\,-12.5,*,-12,\,-11.5$
used in Figure \ref{fig:pp_WLL_pattern} all lie below the instability threshold and belong to the diffusion--chemotaxis instability region. Conversely, $\beta=0$ lies in the stable region, confirming that the corresponding reaction--diffusion system does not exhibit diffusion-driven Turing instability.

The importance of this result goes beyond the mere observation of pattern formation. For the parameter values considered here, the homogeneous equilibrium is not only asymptotically stable, but also non-reactive. Therefore, neither reactivity nor classical Turing instability can explain the emergence of the observed spatial structures. The patterns displayed in Figure~\ref{fig:pp_WLL_pattern} are generated exclusively by the chemotactic mechanism, providing a direct numerical illustration of the theoretical result established in Section~\ref{sec:notnecessary}: reactivity of the homogeneous equilibrium is not a necessary condition for chemotaxis-driven pattern formation.

Besides triggering instability, the chemotactic sensitivity also affects both the amplitude and the morphology of the emerging patterns. Close to the stability threshold ($\beta=-11.5$ and $\beta=-12$), the solution organizes into a nearly regular hexagonal arrangement with relatively small amplitude. As $\beta$ decreases and the system moves deeper into the instability region, the amplitude increases and the hexagonal pattern progressively loses its regularity. For larger negative values of $\beta$ ($\beta=-13.5$ and $\beta=-15$), elongated stripe-like structures become dominant, revealing a morphology transition from hexagonal spot-like patterns to stripes.
This behavior closely resembles the scenario previously observed for the chemotaxis-extended Schnakenberg model. In both examples, moving away from the stability boundary leads to an increase in the pattern amplitude and to a progressive reorganization of the spatial structure. The fundamental difference is that, in the present case, the entire pattern-forming regime lies outside both the reactive and the Turing-instability regions, showing that chemotaxis alone can provide an independent route to spatial self-organization.

\section{Conclusions}
\label{sec:conclusions}

In this work we investigated the relationship between reactivity and pattern formation in reaction--diffusion systems with chemotaxis. The starting point was the classical result of Neubert, Caswell and Murray \cite{Neubert2002}, stating that reactivity of the spatially homogeneous equilibrium is a necessary condition for diffusion-driven, or Turing, instability. We showed that this conclusion does not extend in a direct way to reaction--diffusion--chemotaxis systems.

For a general two-species model, we recalled the linear instability conditions in the presence of both diffusion and chemotactic transport. The analysis shows that chemotaxis modifies the instability mechanism through a non-symmetric contribution to the linearized operator. As a consequence, the necessary condition for instability is no longer determined by the reactivity of the kinetic Jacobian $J^*$, but involves the spatially dependent operator containing the chemotactic contribution. This leads to the main theoretical conclusion of the paper: while reactivity of the homogeneous equilibrium remains necessary for classical Turing instability, it is not necessary for chemotaxis-driven instability.

The geometric interpretation of the instability condition further clarifies this distinction. In the purely diffusive case, a stable non-reactive equilibrium cannot be destabilized by diffusion alone. In contrast, when chemotaxis is present, sufficiently strong directed transport may move the system across the instability boundary even if the local kinetics remains stable and non-reactive. Thus, chemotaxis should not be interpreted merely as a perturbation of the classical Turing mechanism. Rather, it provides an independent route to spatial self-organization, driven by the non-symmetric structure of the taxis operator.

From the computational point of view, we extended the matrix-oriented approach for reaction--diffusion systems, proposed in \cite{DautiliaSguraSimoncini2020}, to the reaction--diffusion--chemotaxis setting. The chemotactic flux was discretized in a form compatible with the ECDF approximation of the Laplacian \cite{SGURA20124132}, so that in the constant-mobility limit the discrete chemotactic operator reduces to the same matrix Laplacian used for diffusion. This construction preserves the Sylvester-equation structure of the implicit diffusion step. Therefore, the inclusion of chemotaxis only modifies the explicit right-hand side of the IMEX scheme, while retaining the computational advantages of the matrix-oriented formulation.

The numerical experiments illustrated two complementary roles of chemotaxis. In the chemotaxis-extended Schnakenberg model, the underlying reaction--diffusion system already exhibits classical Turing instability. In this case, chemotaxis modifies pre-existing diffusion-driven patterns: depending on the sign and magnitude of the chemotactic sensitivity, it may enhance, weaken, suppress, or reshape the spatial structures. The simulations showed that chemotaxis affects both the amplitude and the morphology of the patterns, producing transitions between spots, stripes, and mixed configurations.

The predator--prey example addressed the main theoretical message of the paper. For the selected parameter regime, the coexistence equilibrium is asymptotically stable and non-reactive, as confirmed by the negativity of both the spectral abscissa $\alpha(J^*)$ and the numerical abscissa $r(J^*)$. Moreover, the corresponding reaction--diffusion system does not satisfy the classical Turing instability condition. Nevertheless, when the chemotactic sensitivity is sufficiently negative, spatially heterogeneous stationary solutions emerge. These patterns are therefore generated solely by the chemotactic contribution. This provides a direct numerical confirmation that reactivity of the homogeneous equilibrium is not a necessary condition for chemotaxis-driven pattern formation.

The predator--prey simulations also showed that chemotaxis controls not only the onset of instability, but also the morphology and amplitude of the resulting patterns. Close to the chemotactic stability threshold, low-amplitude hexagonal structures are selected. Moving deeper into the instability region, the amplitude increases and the hexagonal arrangement progressively gives way to stripe-like patterns. This behavior resembles the morphology transitions observed in the Schnakenberg example, but with a key difference: in the predator--prey case, the entire pattern-forming mechanism is induced by chemotaxis in a regime where neither local reactivity nor diffusion-driven instability is present.

Overall, our results highlight a qualitative difference between diffusion-driven and chemotaxis-driven mechanisms of pattern formation. In reaction--diffusion systems, instability necessarily originates from a reactive homogeneous equilibrium. In reaction--diffusion--chemotaxis systems, instead, non-symmetric directed transport can destabilize stable non-reactive equilibria and generate spatial organization through a purely chemotactic mechanism. This suggests that reactivity-based criteria developed for classical Turing systems must be reconsidered when taxis or more general non-symmetric transport processes are present.

Several directions remain open for future work. From the numerical point of view, a first natural extension concerns the development of higher-order time discretizations for reaction--diffusion--chemotaxis systems. In the present work, we have employed a first-order IMEX Euler scheme in order to focus on the construction of an ECDF-compatible chemotactic discretization and on its integration within the matrix-oriented framework. However, long-time simulations of pattern-forming systems may benefit from higher-order stabilized methods. In this direction, Runge--Kutta TASE methods for reaction--diffusion problems, such as those developed by Conte et al.~\cite{ConteMontijanoPaganoPaternosterRandez2026}, represent a promising class of schemes for improving temporal accuracy while retaining favorable stability properties in stiff regimes.

Another relevant direction concerns the qualitative fidelity of the fully discrete method with respect to the instability mechanism. In the reaction--diffusion setting, the preservation of the continuous Turing region by fully discrete matrix-oriented splitting schemes has been recently investigated in~\cite{Monti2026discrete}. There, it was shown that time discretization may generate spurious numerical patterns or, conversely, suppress continuous instabilities. The present results suggest an analogous question for reaction--diffusion--chemotaxis systems: whether a fully discrete scheme preserves the diffusion--chemotaxis instability region and correctly distinguishes chemotactic patterns predicted by continuous theory from numerical artefacts.

A further computational direction is the combination of the matrix-oriented framework with model order reduction techniques. Although matrix-oriented formulations are already significantly more efficient than fully vectorized approaches, large spatial domains, fine grids, and long-time integrations may still lead to substantial computational costs. Reduced-order strategies could therefore be used to further decrease the online cost while retaining the ability to capture pattern-forming dynamics. In the vector-oriented setting, physics-based and data-driven reduction strategies have already been explored for reaction--diffusion pattern formation in~\cite{BMS2021,AMS2023,AMS2024}. A first step toward a matrix-oriented reduced formulation is provided by the two-sided POD--DEIM approach of Kirsten and Simoncini~\cite{KirstenSimoncini2021}, which preserves the matrix structure of semilinear matrix differential equations. Extending these ideas to reaction--diffusion--chemotaxis systems would be particularly useful for simulating chemotaxis-driven patterns on larger domains and increasingly refined grids at reduced computational cost.

From the modelling and analytical point of view, it would be interesting to extend the present analysis to systems with nonlinear taxis sensitivities, cross-diffusion effects, or more than two interacting species. Such extensions would help clarify whether chemotaxis-driven pattern formation from stable non-reactive equilibria is a robust feature of non-symmetric transport processes in biological and ecological systems.

\section*{Acknowledgements}

\noindent A. M. research activity is funded by PR PUGLIA FESR FSE+ 2021--2027 - Fondo Europeo Sviluppo Regionale - Asse Prioritario I ``Competitività e Innovazione'' - Obiettivo specifico RSO1.1 - Azione 1.5 ``Interventi per il rafforzamento del sistema innovativo regionale e sostegno alla collaborazione tra imprese e strutture di ricerca'' - Sub-Azione 1.5.1 ``Supporto alle attività di ricerca e sviluppo su aree tematiche di rilievo e all'applicazione di soluzioni tecnologiche funzionali alla realizzazione delle strategie di S3'', Avviso pubblico ``Reti - Sostegno alla ricerca collaborativa'' approvato con A.D. n. 208/2024, A.D. n. 216/2024, A.D. n. 227/2024, A.D. n. 230/2024 e AD n.3/2025, CUP B89J24003660007, project title ``PRISM - Diagnosi precoce e responsiva nei Vigneti Mediterranei''. 

A. M. research activity is funded by CNR, Project ID: DIT.AD021.228, CUP: B97G25000730001, project title ``MODAMB - Modelli e 
metodi matematici per l'ambiente''.

A. M. is member of the INdAM research group GNCS.

\bibliographystyle{siam}
\bibliography{biblio}

\end{document}